\documentclass[11pt,leqno]{article}
\usepackage{amsthm,amsfonts,amssymb,epsfig,graphics,amsmath,eufrak}
 \usepackage[latin1]{inputenc}\relax

\newcommand{\trans}{\text{\rm tr}}

\newcommand{\sgn}{\text{\rm sgn}}

\newcommand{\RR}{{\mathbb R}}

\newcommand\br{\begin{rem}}
\newcommand\er{\end{rem}}
\newcommand\bp{\begin{pmatrix}}
\newcommand\ep{\end{pmatrix}}
\newcommand\be{\begin{equation}}
\newcommand\ee{\end{equation}}
\newcommand\ba{\begin{equation}\begin{aligned}}
\newcommand\ea{\end{aligned}\end{equation}}

\newcommand\adots{\mathinner{\mkern2mu\raise1pt\hbox{.}
\mkern3mu\raise4pt\hbox{.}\mkern1mu\raise7pt\hbox{.}}}

\newtheorem{theo}{Theorem}[section]

\newtheorem{cor}[theo]{Corollary}
\newtheorem{lem}[theo]{Lemma}
\newtheorem{defi}[theo]{Definition}

\newtheorem{rem}[theo]{Remark}

\numberwithin{equation}{section}
\title{
A sharp stability criterion
for soliton-type propagating phase boundaries in Korteweg's model
}

\author{Kevin Zumbrun\footnote{Indiana University, partially 
supported  by NSF grant DMS-0300487.
Thanks to Bj\"orn Sandstede for pointing out
related results in \cite{BD1, BD2}, and to Thomas Bridges
and Gianne Derks for helpful discussions 
on the relations between their multi-symplectic approach 
and the ``direct method'' followed here.
 }}

\date{May 2, 2006}

\begin{document}

\maketitle
\begin{abstract}
Recently, Benzoni--Gavage, Danchin, Descombes, and Jamet
have given a sufficient condition for linear and nonlinear
stability of solitary wave solutions of Korteweg's 
model for phase-transitional isentropic gas dynamics 
in terms of convexity of a certain ``moment of instability'' with
respect to wave speed, which is equivalent to variational stability
with respect to the associated Hamiltonian energy
under a partial subset of the constraints of motion; they conjecture
that this condition is also necessary.
Here, we compute a sharp criterion for spectral stability in terms of
the second derivative of the Evans function at the origin, and
show that it is equivalent to the variational condition
obtained by Benzoni--Gavage et al, answering their conjecture
in the positive.
\end{abstract}

%%%%%%%%%%%%%%%%%
\bigbreak
\section{Introduction}\label{s1}

Motivated by recent work of Benzoni--Gavage et al \cite{BDDJ},
we investigate in this paper 
stability of ``soliton-type''
(i.e., homoclinic) traveling wave solutions
\be\label{1}
U(x,t)=\bar U(x-st),
\quad \lim_{z\to \pm \infty} \bar U(z) = U_\infty,
\ee
$U=(v,u)$, $\bar U=(\bar v, \bar u)$ of
the Korteweg model
\ba\label{2}
v_t - u_x&=0\\
u_t + p(v)_x&=
%-(\kappa v_{xx})_x - (1/2)\kappa_v (v_x)^2,\\
-\kappa v_{xxx},\\
\ea
for isentropic phase-transitional gas dynamics, written in Lagrangian
coordinates, with $v$ denoting specific volume, $u$ particle velocity,
$p$ pressure, and $\kappa>0$ a coefficient of capillarity, taken
for simplicity to be constant.
The extension of our results to general $\kappa(v)$ as considered 
in \cite{BDDJ} is straightforward; see Remarks \ref{rigor} and \ref{signrem}.
As pointed out in \cite{BSS}, system \eqref{2} is also formally equivalent to
the ``good'' Boussinesq equation modeling shallow-water flow.

Equations \eqref{1}, accounting for the effects
of compressibility and capillarity, but neglecting viscosity,
are of dispersive, type, in contrast to the dissipative type
of the usual compressible Navier--Stokes equations. 
Indeed, as discussed in greater generality in \cite{BDDJ},
\eqref{1} has the Hamiltonian structure
\be\label{ham}
U_t= {\cal J} \delta \cal{H},
\ee
where 
\be\label{J}
{\cal J}=\partial_x J:=\partial_x \bp 0 & 1\\ 1 & 0 \ep
\ee
 is a constant-coefficient skew-symmetric first-order differential 
operator and $\delta \cal{H}$ is a second-order differential operator 
corresponding to the variational derivative of the Hamiltonian functional
\be\label{H}
{\cal H}=\int H,
\quad H=(1/2)(u-u_\infty)^2 -\int_{v_\infty}^v (p(z)-p(v_\infty))dz 
+ (1/2)\kappa(v) (v_x)^2,
\ee
of the (relative) total energy $H$ of the system, with
$\kappa(v)\equiv \kappa$ constant.
%Same form also for general $\kappa(v)$...

Formally, 
$$
(d/dt){\cal H}(U)=
\langle \delta {\cal H}, U_t\rangle=
\langle \delta {\cal H}, {\cal J} \delta {\cal H}\rangle=0,
$$
so that the Hamiltonian is one conserved quantity of motion.
A second (formally) conserved quantity, arising as a consequence of group
invariance under translation (see \cite{GSS}) is the relative generalized
momentum
\be\label{Q}
{\cal Q}(U)= (1/2) \langle J(U-U_\infty), (U-U_\infty)\rangle=
\int (u-u_\infty)(v-v_\infty)(x) dx
\ee
(formally, $J=(\partial_x J)^{-1}\partial_x$, as prescribed in \cite{GSS}).
Two additional conserved quantities
are the relative masses 
$$
{\cal P}_1(U)=\int (v-v_\infty)(x)dx, 
\quad
{\cal P}_1(U)=\int (u-u_\infty)(x)dx.
$$
The existence of these further quantities is
associated with the fact that operator $\cal J$ is not onto,
a circumstance that turns out to be significant. 

With this framework, it was shown by Benzoni--Gavage et al \cite{BDDJ} that
stability of solitons may be investigated
by variational methods, following the formalism of \cite{GSS, BS}.
Specifically, 
one may compute that solitary wave solutions are
critical points of the Hamiltonian ${\cal H}$ under constraint
${\cal Q}$, satisfying the Euler--Lagrange equation (itself Hamiltonian)
\be\label{EL}
(\delta {\cal H} + s \delta {\cal Q})
(\bar U)=0,
\ee
where 
%$\nu$ is a Lagrange multiplier equal to the speed $s$.
speed $s$ plays the role of a Lagrange multiplier.
Such solutions occur in a one-parameter family $\bar U^s$.
Formally, {\it strict variational stability} of $\bar U$ with 
respect to constraint ${\cal Q}$
is thus sufficient for {\it time-evolutionary orbital stability}
of the family $\{\bar U^s\}$, since, then, (i) A minimum $\bar U^{s({\cal Q})}$
should therefore persist under small changes in ${\cal Q}$, and
(ii) within each level surface of ${\cal Q}$, 
the conserved quantity
${\cal H}(U)-{\cal H}(\bar U^{s({\cal Q})})$ 
should control $\|U-\bar U^{s(\cal Q})\|^2$ in
the underlying Hilbert norm.

In \cite{GSS}, Grillakis et al made rigorous this intuitive argument, showing
that variational stability indeed implies nonlinear orbital stability in
a quite general framework (which applies here).  Moreover, 
under the assumption
that ${\cal L}^s$ {\it have at most one negative eigenvalue}, 
they gave
a necessary and sufficient condition for strict variational stability
in terms of strict convexity with respect to $s$ 
of the ``moment of instability'' 
\be\label{moment}
m(s):= ({\cal H}+ s {\cal Q})(\bar U^s)
\ee
(in the terminology of \cite{BS}),
i.e.,
%\be
%\label{convex}
%d^2m/ds^2= -\langle \partial \bar U^s/\partial s, {\cal L}^s
%\partial \bar U^s/\partial s\rangle> 0,
%\ee
\be\label{monotone}
(d^2m/ds^2)(\bar s)=
(d/ds) {\cal Q}(\bar U^s)=
\langle (\bar U-U_\infty), J(\partial \bar U^{\bar s}/\partial s)\rangle>0.
\ee
%, under the assumption
%that ${\cal L}^s$ has {\it at most one negative eigenvalue}. 
Finally, for $\cal J$ {\it onto}, they showed by construction of
a suitable Lyapunov function that strict failure of
\eqref{convex} implies time-evolutionary instability, completely
deciding the issue of stability.

\br\label{mon}
\textup{
%so that \eqref{convex} is equivalent to
%$(d/ds) {\cal Q}(\bar U^s)>0$;
Relation \eqref{monotone} shows that
convexity of the moment of instability is equivalent
to monotone increase with respect to $s$
of the generalized momentum $\cal Q$ along $\{\bar U^s\}$.\footnote{
In the notation of \cite{GSS}, monotone decrease with respect to $\omega=-s$.}
Local monotonicity of $\cal Q$ along $\{\bar U^s\}$ is 
%a necessary condition for orbital stability, since the map from
necessary for the picture of stability described above
(in particular, for strict variational stability),
since the map from
$s$ to ${\cal Q}(\bar U^s)$ must be locally invertible if all
sufficiently small perturbations
(corresponding to small variations in $\cal Q$) 
are to lie near some $\bar U^s$. 
Through the key relation 
\be\label{role}
{\cal L}^s (\partial \bar U^s/\partial s)= -\delta Q
\ee
(obtained by differentiating \eqref{EL} with respect to $s$;
see Proposition 4, \cite{BDDJ}), we find that \eqref{monotone}
is equivalent to
\be
\label{convex}
d^2m/ds^2= -\langle \partial \bar U^s/\partial s, {\cal L}^s
\partial \bar U^s/\partial s\rangle,
%> 0,
\ee
where the self-adjoint operator ${\cal L}^s$, defined by
\be\label{hess}
{\cal L}^s=(\delta^2 {\cal H}+ s \delta^2 {\cal Q})|_{\bar U^s},
\ee
is the constrained Hessian about $\bar U^s$.
}
\er

In \cite{BDDJ}, a simple formula is given for $m(s)$ and 
evaluated numerically to show that 
regions of both convexity and nonconvexity of $m(\cdot)$
may arise, depending on physical parameters.
The first case corresponds with orbital stability, as discussed
above.  However, {\it since $\cal J$ is not onto}, the second
case is inconclusive by the theory of \cite{GSS}. 
The authors conjecture nonetheless that convexity is necessary as well as 
sufficient for stability, so that the second case in fact 
corresponds to instability.  

Here, we investigate time-evolutionary stability directly, using
alternative, Evans function methods introduced in
\cite{AGJ, PW, GZ, Z, BSZ} to obtain a sharp criterion
for spectral time-evolutionary stability in terms of the sign
of the second derivative at the origin of the Evans function $D(\lambda)$,
an analytic function whose zeroes correspond in location and multiplicity
with eigenvalues of $L$.
Specifically,
%denoting the Evans function by $D(\lambda)$, 
we obtain the following results 
deciding in the positive the conjecture of \cite{BDDJ}.

\begin{theo}\label{main}
$D''(0)= c(d^2m/ds^2)(\bar s)$ for a nonzero constant $c$.
\end{theo}

\begin{cor}\label{stab}
Traveling waves \eqref{1} of \eqref{2} are 
time-evolutionarily (linearly and nonlinearly)
stable if $d^2m/ds^2> 0$ and only if $d^2m/ds^2\ge 0$,
i.e., if they are strictly variationally stable with respect to
constraint $\cal Q$,
and only if they are nonstrictly variationally stable.
\end{cor}

Theorem \ref{main} generalizes a number of similar results
obtained in \cite{PW} for various related scalar models, and
indeed holds in far greater generality for systems of abstract
form \eqref{ham}.
To see why, and to better understand in general the relation between
the moment condition and
variational and time-evolutionary stability, notice that,
in this generality, \eqref{Q} and\eqref{monotone} become, formally
(see \cite{GSS}), 
$$
{\cal Q}=(1/2)\langle {\cal J}^{-1} \partial_x (U-U_\infty),
(U-U_\infty \rangle
$$
and 
$$
(d^2m/ds^2)(\bar s)= 
\langle {\cal J}^{-1}\partial_x (\bar U-U_\infty), (\partial \bar U^{\bar s}/\partial s)\rangle 
=
\langle {\cal J}^{-1}\partial_x \bar U, (\partial \bar U^{\bar s}/\partial s)\rangle. 
$$
On the other hand, differentiation of the traveling-wave equation with
respect to $x$ and $s$, respectively, reveals that
$f_1=\partial_x \bar U$ 
is a right zero-eigenfunction of the linearized operator $L$ and 
$f_2=-\partial_s \bar U^{\bar s}$ a generalized zero-eigenfunction of
height two; for further discussion, see Section \ref{prelim}.
Noting that left and right zero eigenfunctions $\tilde f_1$ and $f_1$ 
of $L=\cal J L$ are related, formally, by $\tilde f_1={\cal J}^{-1}f_1$,
we find the general relation
\be\label{toddrel}
(d^2m/ds^2)(\bar s)= 
-\langle \tilde f_1, f_2 \rangle.
\ee

But, vanishing of the inner product \eqref{toddrel} 
of the genuine left eigenfunction $\tilde f_1$
against the generalized right eigenfunction $f_1$ of height two 
precisely detects existence of a generalized eigenfunction of height three
(by Jordan form), i.e., algebraic multiplicity of order two or more
of the eigenvalue $\lambda=0$.
Thus, 
$(d^2m/ds^2)(\bar s)$
must be a nonzero multiple of $D''(0)$, since vanishing
of $D''(0)$ detects the same phenomenon.

Formula \eqref{toddrel}
%, to our knowledge new in the literature,
reveals a direct (formal) link between the moment condition 
and time-evolutionary stability 
that, moreover, does not require
invertibility of $\cal J$ on the whole space, but only on the range
of $\partial_x$,
%NOTES: Just need ${\cal J}^{-1}\partial_x$ defined...
through the distinguished roles of 
$(\partial \bar U^{\bar s}/\partial x)$ and 
$-(\partial \bar U^{\bar s}/\partial s)$ as zero eigenfunction and
generalized eigenfunction of $L$.
Recall that these same functions played critical roles
in the argument of \cite{GSS} 
linking the moment condition and variational stability:
$(\partial \bar U^{\bar s}/\partial x)$ as zero eigenfuction of $\cal L$
and
$-(\partial \bar U^{\bar s}/\partial s)$ through relation \eqref{role}.

This formal argument can be made rigorous using the extended
spectral theory of \cite{ZH}, Section 6,
valid essentially wherever an Evans function can be analytically 
defined,\footnote{
See also \cite{MZ, Z2} for extensions to operators $L$ of
the degenerate type considered here.
}
concerning Jordan structure at eigenvalues embedded in essential spectrum,
along with additional assumptions assuring that there exist
no additional (extended) genuine zero eigenfunctions other
than $f_0$; see Remark \ref{rigor}.
We shall instead follow the more concrete approach of direct
Evans function calculations using the specific structure at hand,
which provide at the same time sign information.
However, we note that
these are quite similar to those on which
the abstract development of \cite{ZH} is based
(see in particular the proof of Proposition 6.3, \cite{ZH}).

\begin{rem}\label{kapitula}
\textup{
Along similar lines, the general relation
\be\label{kform}
\partial_\lambda^k D(0)=c \langle \tilde f_1, f_{k}\rangle
\ee
has been established by Kapitula
\cite{K} for general (not necessarily Hamiltonian) systems for
which $\lambda=0$ is an isolated eigenvalue of geometric multiplicity
one of the linearized operator $L$ under consideration.
In the present (Hamiltonian) case, $\lambda=0$ is embedded in
the essential spectrum of $L$, and so this result does not apply.
Indeed, Benzoni--Gavage, Serre, and Zumbrun \cite{BSZ} have shown
in the context of viscous conservation laws
that \eqref{kform} does not hold in general for embedded eigenvalue
$\lambda=0$, but rather must be corrected by the addition of 
appropriate boundary terms at plus and minus spatial infinity;
see also related results in \cite{KS1, KS2, KR} for perturbed NLS equations.
Thus, the argument above reflects partly 
the special features of the Hamiltonian case.
}
\end{rem}

\begin{rem}\label{rigor}
\textup{When ${\cal J}$ is a differential operator 
of the general form $\partial_x J$ considered
here, with $J$ an invertible (not necessarily differential) operator 
onto $L^2$ (the most general case falling under the sufficient but
not necessary theory of \cite{GSS}),
we may take ${\cal J}^{-1}\partial_x=J^{-1}$ in the discussion above, and 
$$
\tilde f_1={\cal J}^{-1}f_1
= {\cal J}^{-1}\partial_x (\bar U- U_\infty) 
= J^{-1} (\bar U-U_\infty).
$$
Noting that $\tilde f_1$ decays exponentially as $x\to \pm \infty$,
we find that it is indeed a left genuine zero-eigenfunction of $L$,
making rigorous sense of relation \eqref{toddrel}.
%}
%
%\textup{
To complete the rest of the formal argument sketched above, 
note that bounded solutions of $Lf=0$
correspond to nearby homoclinic connections with different endstates,
hence have equal limits as $x\to \pm \infty$, and this eliminates
them as possible extended eigenfunctions (which might in general
be merely bounded \cite{ZH}) unless they in fact vanish at $\pm \infty$.
Thus, the standard assumption that $\bar U$ be a transverse connection,
ensuring that the $L^2$ kernel of $L$ is one-dimensional,
is sufficient to ensure that the extended kernel of $L$ is also one-dimensional,
i.e., there is a unique (up to constant factor) extended genuine eigenfunction
$f_1$ with dual eigendirection $\tilde f_1$.
}

\textup{
Moreover, since 
$\tilde f_1$, $f_1$, and $f_2$ are all exponentially decaying,
there exists an extended generalized eigenfunction $f_3$, bounded
but not necessarily decaying at infinity, if and only if
$\langle \tilde f_1, f_2\rangle=0$.  For, $Lf_3=f_2$ implies
$$
\langle \tilde f_1, f_2\rangle=
\langle \tilde f_1, Lf_3\rangle= \langle L^*f_1, f_3\rangle=0.
$$
On the other hand, if the extended Jordan chain is order two,
then (see \cite{ZH}, Prop. 5.3 (iii)) 
the extended spectral projection 
$$
{\cal P}=
f_1 \langle \tilde f_2, \cdot\rangle 
+ f_2 \langle \tilde f_1, \cdot\rangle
$$
preserves exponentially decaying elements of the extended eigenspace; in
particular ${\cal P}f_2=f_2$, whence $\langle \tilde f_1, f_2\rangle \ne 0$.
Finally, recalling (\cite{ZH}, Theorem 6.3 (ii)) that the Evans function
vanishes with multiplicity equal to the dimension of the corresponding
extended eigenspace, we obtain the result 
\be\label{weak}
%$D''(0)= c(d^2m/ds^2)(\bar s)$, $c\ne0$
D''(0)= c(d^2m/ds^2)(\bar s), \quad c\ne0
\ee
for this general class: in particular, for the general 
isentropic Korteweg equations discussed in \cite{BDDJ} consisting of
\eqref{ham}--\eqref{H} with arbitrary $\kappa(v)>0$.
%TODO: put back somewhere????
%It is worth noting that conclusion \eqref{weak}
%holds {\it independent of any assumptions on the spectrum of ${\cal L}^s$},
%so in this sense 
%generalizes the results of \cite{GSS}.
}
\end{rem}

\begin{rem}\label{dual}
\textup{
The same extended spectrum argument used in Remark \ref{rigor},
applied to general (not necessarily Hamiltonian) PDE,
yields the remarkable fact that the dual version
\be\label{dualkform}
\partial_\lambda^k D(0)=\tilde c \langle \tilde f_k, f_1 \rangle,
\ee
$\tilde c\ne 0$, of \eqref{kform} remains valid 
in complete generality, {\it without additional boundary terms},
for embedded eigenvalues with (extended) geometric multiplicity
one for which the associated eigenfunction $f_1$ is exponentially
decaying, a condition that is in the traveling wave context 
essentially always satisfied.\footnote{Likewise, the obvious
extensions to higher geometric multiplicity remain valid
so long as all genuine extended eigenfunctions 
are exponentially decaying.}
%NOTE: product of similar inner products, provided chains of
%different length.  If chains of same length, get gramian determinant.
For, then ${\cal P}f_1=f_1$, where
${\cal P}= \sum_{j=1}^K f_j \langle \tilde f_{K-j+1}, \cdot\rangle$
and $K\ge k$ is the order of vanishing of the Evans function, whence
$\langle \tilde f_k, f_1\rangle=0$ if and only if $K \ge k+1$, or,
equivalently, $\partial_\lambda^k D(0)=0$.
Boundary terms arise in the forward formula through the integration
by parts converting
$\langle \tilde f_k, f_1 \rangle=  \langle \tilde f_k, L^{k-1} f_k \rangle$ 
to
$-1^{k-1}\langle (L^*)^{k-1}\tilde f_k, f_k \rangle= 
-1^{k-1}\langle \tilde f_1, f_k \rangle$, except in the case, as in Remark
\ref{rigor} above,
that $f_1, \dots, f_k$ or $\tilde f_1, \dots, \tilde f_k$ 
decay exponentially.  
(Recall that extended eigenfunctions $f_j$, $\tilde f_j$
do not necessarily decay at infinity, except, by assumption, $f_1$,
but rather grow at most algebraically \cite{ZH, BSZ}.)
}
\end{rem}

\begin{rem}\label{signrem}
\textup{
The relation \eqref{weak} by itself does not give
give the complete stability information of Corollary \ref{stab}, 
but requires further information on the {\it sign of $c$} 
(with $D$ suitably normalized: more precisely, the
relation between $\sgn (c)$ and $\sgn D(+\infty)$).
In practice, this is often not restrictive,
since we may calibrate the sign of $c$ 
by continuation (homotopy) to a case where stability is
decided by the theory of \cite{GSS}, i.e., either 
$d^2m/ds^2>0$ (variational stability) or $\cal J$ onto.
For example, in the case considered here, the collection of
all soliton solutions of all isentropic Korteweg models
comprises an open set in parameter space, so we may conclude
by the (numerically observed) existence of stable waves.
Alternatively, the regular perturbation
${\cal J}^\theta:= \partial_x J + \theta K$,
$K:=\bp 0 & 1 \\ -1 & 0 \ep$,
$\theta $ sufficiently small,
preserves sign information while converting $\cal J$
to an operator ${\cal J}^\theta$ that is onto for all $\theta\ne0$.
}

\textup{
A still more general approach, not limited to the Hamiltonian setting,
is to work directly from \eqref{dualkform}, computing the sign
by a direct calculation like that of Section \ref{s:evans}.
In practice, one may sometimes determine the normalization
without doing the full calculation,
in which case it is easier to determine $\tilde f_k$, $f_1$ numerically;
see, for example, the analysis of stability of undercompressive
traveling waves of thin-film models in \cite{BMSZ} (Prop. 2.11
and footnote 6 on duality).
We point out that the righthand side of \eqref{dualkform}
may by itself be considered as a generalized Melnikov integral,
like $d^2m/ds^2$ giving geometric information about the
dynamics of the traveling-wave ODE; see
\cite{MZ} Section 4.2 (especially eqns. (4.9)--(4.13))
for a general duality principle linking
dual eigenfunctions to solutions of the adjoint ODE.  
}
\end{rem}

\begin{rem}\label{BDDJ}
\textup{
There appears in \cite{BDDJ} the statement, 
apparently contradicting Corollary \ref{stab}, 
that there exist stable solitons that are variationally unstable.  
However, the time-evolutionarily stable waves considered in
\cite{BDDJ} are in the standard sense variationally stable; indeed, this is
the property that is used to prove time-evolutionary stability.
The instability that is referred to, rather, is of the constrained
Hessian, \eqref{hess}, with respect to {\it unconstrained} perturbations,
a notion that is necessary but not sufficient for instability with
respect to constraint $\cal Q$.
}
\end{rem}

{\bf Note.} After the completion of this paper, we have learned
of results of Bridges and Derks \cite{BD1, BD2} establishing the 
relation $D''(0)=c(d^2m/ds^2)$ for the good Boussinesq system,
which is equivalent to the main example \eqref{2} studied here.
More generally, they derive 
%general 
formulae generalizing those of
Pego and Weinstein \cite{PW} for scalar equations for the first
nonvanishing derivative of a ``symplectic Evans function'' at the origin,
for systems that can be put in the multi-symplectic form
$MZ_t+KZ_x=\nabla S(Z)$, where $Z\in{\Bbb R}^{2n}$, $M$, $K$ 
are constant skew-symmetric $2n\times 2n$ matrices, and $S$ is a 
smooth function on ${\Bbb R}^{2n}$.
These include in principle a rather large class of Hamiltonian PDE,
overlapping with but apparently distinct from the class 
(described in Remark \ref{rigor}) to which our methods apply.
%TODO: old, delete?
%Note 
%(Remark \ref{rigor}) 
%that our methods apply 
%to the general (arbitrary $\kappa(v)$) isentropic 
%equations studied in \cite{BDDJ}, and not only \eqref{2};
%in addition, they are considerably more straightforward than those of 
%\cite{BD1, BD2},
%being carried out in the original coordinates and motivated by the
%simple relation \eqref{toddrel}.

The multi-symplectic approach has the advantage that
it is equation-independent once a change of coordinates
to multi-symplectic form has been found, yielding automatically
a characterization of the sign of $c$ in \eqref{weak} in terms
of the geometry of the phase space of the traveling-wave ODE
(cf.  Remark \ref{signrem}, par. two).
%.\footnote{
%Equivalent, of course, to the 
%one that is ultimately generated by the equation-specific
%calculations described in Remark \ref{signrem}, paragraph two.}
On the other hand, our methods are somewhat more straightforward,
being carried out in the original coordinates and motivated by the
simple relation \eqref{toddrel}.
In addition (see Remarks \ref{dual}--\ref{signrem}),
they yield useful partial information also in the general, 
non-Hamiltonian case.
%the underlying ideas extend also to the non-Hamiltonian case.
%TODO: delete, I think- numerically could be quite easily done... so,
%indeed IS useful...
%The characterization of $\sgn(c)$ by either method is largely theoretical,
%since it cannot be computed analytically except in very simple cases.
%TODO: old, delete?
%In practice, determination of $\sgn(c)$ is intractable by either
%approach except in very simple cases.
%NOTE: specifically, consists in dot product of $\partial \bar U^s/\partial s$
%against a special left e-vector selected by Evans calculation/phase space
%geometry as dual to collection of right decaying solutions of e-value
%equation.  In general, unknowable even if one derives it except in simple
%cases or as above by comparison to GSS...

\section{Preliminaries}\label{prelim}

Substituting \eqref{1} into \eqref{2}, we obtain the traveling-wave equation
\ba\label{trav}
-sv' - u'&=0\\
-su' + p(v)'&= -\kappa v''',\\
\ea
or, substituting the first equation into the second,
and integrating from $-\infty$, the Hamiltonian ODE
(nonlinear oscillator)
\be\label{inttrav}
v''=
\kappa^{-1}\big(s^2v + p(v)- s^2v_\infty - p(v_\infty)\big).
\ee
Alternatively, we may write \eqref{trav} formally as
$-sU'={\cal J}{\cal H}(U)$, or
\be\label{alttrav}
{\cal J}(\delta {\cal H}+sJ)(U)=
{\cal J}\delta( {\cal H}+s{\cal Q})(U)= 0.
\ee

By the Hamiltonian structure of \eqref{inttrav}, it follows
that homoclinic orbits persist under changes in speed $s$ and endstate
$U_\infty$, forming for fixed $U_\infty$
a one-parameter family $\bar U^s$, $s\ge 0$, as described in the
introduction.
The equations are invariant under shifts in velocity $u$, so that
any value of $u_\infty $ is possible.  However, the requirement
that $U_\infty$ be a saddle-point of \eqref{trav} enforces
on $v_\infty$ the condition
\be\label{signalpha}
s^2+ p'(v_{\infty})< 0.
\ee

Making the standard change of coordinates $x\to x-st$ 
to a rest frame for the traveling wave, we may investigate
its stability as an equilibrium solution 
$U(x,t)=\Bar U(x)$ of
$U_t-sU_x={\cal J}\delta{\cal H}(U)$, or
\be\label{nonlin}
U_t={\cal J}\delta({\cal H}+ s{\cal Q})(U).
\ee
Linearizing \eqref{nonlin} about $\bar U$, we obtain
\be\label{lin}
U_t=LU:={\cal J}{\cal L}^{\bar s},
\ee
where ${\cal L}^s$ is defined as in \eqref{hess} and $\bar s$
denotes the speed of the wave $\bar U$ under investigation.

\begin{defi} \label{stable}
\textup{
The wave $\bar U$ is spectrally stable if the spectrum
$\sigma(L)$ of the linearized operator $L$ about the wave
is contained in the nonpositive complex half-plane
$\{\lambda:\, \Re \lambda \le 0\}$.
}
\end{defi}

Routine calculation (see, e.g., \cite{BDDJ, GSS, PW}, or Section \ref{s:evans}
below)
shows that the essential spectrum of $L$ consists of the entire imaginary axis,
so that stability is at best of neutral, or bounded type, rather than
asymptotic stability.  This may be seen, likewise, by the fact that the
equations are time-space reversible.
Spectral stability is therefore determined by the point spectrum
of $L$: specifically, whether there lie eigenvalues off of the imaginary
axis.
The following general result of Pego and Weinstein,
a quantitative (linear) version of the previously-remarked
relation between variational and time-evolutionary stability,
gives a way to bound the number of such eigenvalues.

\begin{lem} [\cite{PW}] \label{PWbound}
For a linear operator $L$ factoring as $L={\cal J}{\cal L}$
with $\cal J$ skew-symmetric and $\cal L$ self-adjoint,
the number of eigenvalues of $L$ in the positive complex
half-plane $\{\lambda:\, \Re \lambda >0\}$ is less than or
equal to the number of negative eigenvalues of $\cal L$.
\end{lem}

\begin{proof} See \cite{PW}.
\end{proof}

Lemma \ref{PWbound} was used in \cite{BDDJ} to establish the
following upper bound.
This is not needed in the present context, in which we seek to
establish {\it instability}, but gives useful additional information;
see Remark \ref{count}. (It is necessary for stability; see 
Remark \ref{nec}.)

\begin{cor}[\cite{BDDJ}]\label{bound}
The number of unstable (i.e., positive real part) eigenvalues of
$L={\cal J}{\cal L}^{\bar s}$ is less than or equal to one.
\end{cor}

\begin{proof}
Eigenvalues of ${\cal L}^{\bar s}$ may be shown to correspond
to eigenvalues of a (second-order scalar) Sturm--Liouville operator
$Mv:=v''- \kappa^{-1}(s^2 + \alpha)v$, $\alpha:= p'(\bar v(x))$
in the variable $v$.
Since $\bar v'$ by translation invariance is a zero eigenfunction
of this operator 
(see further discussion below; or, just differentiate \eqref{trav}), 
the number of unstable eigenvalues by standard
Sturm--Liouville theory is equal to
the number of nodes (zeroes) of $\bar v'$, which, by \eqref{trav}
may be seen to be one.
See \cite{BDDJ} for details.
\end{proof}

A general approach to resolving the issue of variational vs.
evolutionary stability when variational methods fail to decide
the question, introduced in \cite{PW},
is via the Evans function $D(\lambda)$, an analytic function taking
real values to real values, whose
zeroes correspond to eigenvalues of $L$; for origins of the Evans
function, see \cite{E1, E2, E3, E4, AGJ}.  
By translational invariance,
$D(0)$ necessarily vanishes.  For, differentiating \eqref{alttrav} with
respect to $x$, we obtain $L\bar U'=0$, hence $\bar U'$ is a zero
eigenfunction of $L$.  
Likewise, existence of a one-parameter family $\bar U^s$
of solutions with the same endstate 
implies that $D'(0)=0$.  For, differentiating \eqref{alttrav} with
respect to $s$ yields $L(\partial\bar U^s/\partial s)|_{s=\bar s}=-\bar U'$, 
hence $-(\partial \bar U^s/\partial s)|_{s=\bar s}$ 
is a generalized zero eigenfunction
of $L$ and there is a nontrivial Jordan block for $L$ at $\lambda=0$.

Since $D(+\infty)$ by standard considerations (see \cite{AGJ, PW, GZ}) 
has a constant,
nonzero sign as parameters are varied, this means that, in the generic
case that $D''(0)\ne 0$, the number
of unstable (i.e., positive) real roots of $D$ (eigenvalues of $L$) is odd 
or even depending on the sign of $D''(0)$.  
Used in conjunction with Lemma \ref{PWbound}, this observation can
yield complete information.
Namely, when the total number
is at most one, spectral stability is determined by the sign of
the second derivative of the Evans function.

The second derivative $D''(0)$ was evaluated in \cite{PW}
for scalar KdV-type equations and shown in several related cases to
be exactly the second derivative with respect to $s$ of the moment
of instability, thus establishing necessity of $d^2m/ds^2\ge 0$ 
along with sufficiency of $d^2m/ds^2>0$ for linearized
and nonlinear stability in those cases.
Their method of computing $D$ and its derivatives 
does not apply in the system case considered here.
However, it was shown in \cite{GZ} and \cite{BSZ}, respectively,
that quite similar formulae may be obtained for $D'(0)$ and
higher derivatives $(d/d\lambda)^k D(0)$ in the system case;
see also \cite{KS1, KS2, KR, LP, BD1, BD2}.

In the remainder of the paper, we calculate $D''(0)$ by
the method of \cite{BSZ} and show by explicit computation that 
it is a nonzero multiple of $(d^2m/ds^2)(\bar s)(\bar s)$,
thus establishing Theorem \ref{main} and 
%(eventually) 
Corollary \ref{stab}.

\begin{rem}\label{nec}
\textup{
The relation $(d^2m/ds^2)(\bar s)=cD''(0)$, $c\ne 0$,
shows that the assumption that ${\cal L}^s$ have at most one
negative eigenvalue is necessary for the result of \cite{GSS}
that $(d^2m/ds^2)(\bar s)>0$ implies stability.
For, by the discussion above, the sign of $D''(0)$, hence
of $(d^2m/ds^2)$, counts the parity of the number
of unstable eigenvalues, which might be even in general despite
instability of $L$.
}
\end{rem}

\section{Evans function calculations}\label{s:evans}

\subsection{Construction of the Evans function}
Writing out the eigenvalue equation $(L-\lambda) U=0$, $L$ as in 
\eqref{lin}, we obtain
\ba\label{eval}
\lambda v -sv' - u'&=0\\
\lambda u -su' + (\alpha v)'&= -\kappa v''',\\
\ea
where $\alpha:=p'(\bar v(x))$ as in \eqref{trav}, which
may be written in phase variables $W=(u,v,v',v'')^\trans$ as
a first order system of ODE $W'=A(x,\lambda)W$, or
\be\label{firstorder}
\bp u\\v\\v'\\v''\ep'=
\bp 0 & \lambda & -1 & 0\\
0 & 1 & 0 & 0 \\
0 & 0 & 1 & 0\\
-\lambda & s \lambda & -\alpha -s & -\alpha' \bar v_x
\ep
\bp u\\v\\v'\\v''\ep.
\ee
Using the fact (easily verified for \eqref{trav}) that $\bar U(x)$
converges exponentially to its limit $U_\infty$ as $x\to \pm \infty$,
we may construct an Evans function for \eqref{firstorder} by
the general method described in \cite{GZ}; see \cite{Z2}
for a more recent exposition. 

Examining first the constant-coefficient limit
\ba\label{cc}
\lambda v -sv' - u'&=0\\
\lambda u -su' + \alpha_\infty v'&= -\kappa v'''\\
\ea
of \eqref{eval}, $\alpha_\infty<0$ by \eqref{signalpha},
we find, taking the Fourier transform $\partial_x \to i\xi$, that
the spectrum of limiting, constant coefficient operator $L_\infty$
of $L$ satisfies dispersion relation
\be\label{disp}
\lambda(\xi)=si\xi \pm \sqrt{\alpha_\infty -\kappa \xi^4}, 
\quad \xi\in \RR,
\ee
hence consists of the imaginary axis.
By a standard result of Henry \cite{He} on asymptotically
constant-coefficient operators, we thus find that the essential
spectrum of $L$ likewise consists of the imaginary axis, while
spectra lying in the strictly unstable half-plane 
$\Lambda:=\{\lambda:\, \Re \lambda >0\}$ 
consists entirely of isolated eigenvalues.

The same calculation shows that the dimensions of the stable 
and unstable subspaces of the limiting coefficient matrix 
$A_\infty(\lambda):=\lim_{z\to \pm \infty}A(z, \lambda)$
is {\it constant} in $\Lambda$.
For, substituting $\mu=i\xi$ into the characteristic equation
\be\label{char}
(\lambda-s\mu)^2+\alpha_{-\infty}\mu^2 + \kappa \mu^4=0
\ee
determining eigenvalues $\mu(\lambda)$,
we obtain the dispersion relation \eqref{disp}; thus, eigenvalues
$\mu(\lambda)$ may cross the imaginary axis only along the dispersion
curve, i.e., for $\lambda$ on the imaginary axis.
Taking $\lambda \to +\infty$ along the real axis in \eqref{char}, 
we find that $\lambda^2\sim -\kappa \mu^4$, so that
$\mu$ lie along the fourth roots of $-1$.
Thus, both stable and unstable subspaces of
$A_\infty(\lambda)$ have dimension two in $\Lambda$.
An elementary matrix perturbation calculation at $\lambda=0$ 
shows that these extend analytically to a neighborhood of the origin,
$\lambda=0$, with limiting values at $\lambda=0$ given by the
direct sum of stable (resp. unstable) subspace of $A(\lambda)$
and the vector
\be\label{limvect}
W=(\sqrt{-\alpha_\infty},1,0,0)^\trans
\quad
\hbox{\rm (resp. $(-\sqrt{-\alpha_\infty},1, 0,0)^\trans$},
\ee
with $U=(v,u)$ coordinate equal to the
unstable (resp. stable) subspace 
$$
U=(1, \sqrt{-\alpha_\infty})^\trans
\quad
\hbox{\rm (resp. $(1, -\sqrt{-\alpha_\infty})^\trans$}
$$
of the convection matrix
$$
a:=\bp -s & -1\\
\alpha_\infty & -s\ep.
$$ 
(Recall, \eqref{signalpha}, that $\det a<0$, so that $a$ has
one real positive and one real negative eigenvalue.)

Applying the framework of \cite{GZ}, we find that, on $\Lambda$,
the subspaces of solutions decaying at $+\infty$ and $-\infty$ of the
variable-coefficient equations \eqref{eval} likewise have dimension two.
Moreover, there exist choices of bases $W_1^+, W_2+ $ and
$W_3^-, W_4^-$ for these subspaces that are analytic in $\lambda$
on $\Lambda$ and extend analytically to $\Lambda \cup B(0,r)$
for $r>0$ sufficiently small.
At $\lambda=0$,  
we may choose $W^+_1(0,\cdot)=W_4(0,\cdot)^-=\partial_x \bar W^{\bar s}$ 
asymptotically decaying,  $\bar W^{\bar s}:=
(\bar u^{\bar s},
\bar v^{\bar s},
\bar v^{\bar s}_x,
\bar v^{\bar s}_{xx})^\trans$,
and $W_2^+(0,\cdot)$, $W_3^-(0,\cdot)$ asymptotically constant, 
with
\ba\label{lims}
\lim_{x\to +\infty} W^+_2(x,0)&=( c,1,0,0)^\trans, \\
\lim_{x\to -\infty} W^-_3(x,0)&=
( -c,1,0,0)^\trans,
\ea
$c:= \sqrt{-\alpha_\infty}>0$.
The Evans function is defined as the Wronskian
\be\label{evans}
D(\lambda):=\det \bp W_1^+& W_2^+
& W_3^-& W_4^-\ep|_{x=0},
\ee
zeroes of which detect nontrivial intersection between decaying manifolds 
of solutions at plus/minus spatial infinity, i.e., decaying solutions of 
eigenvalue equation \eqref{eval}.

\subsection{Calculation of $D''(0)$}
%\subsection{Proof of Theorem \ref{main}}

We now compute $D''(0)$ by a simplified version (taking advantage
of special structure) of the approach of \cite{BSZ}.

\begin{proof}[Proof of Theorem \ref{main}]
Following \cite{BSZ} (see also related calculations of \cite{ZH}, Section 6),
we choose a convenient basis for the calculation of derivatives at the origin,
based on the Jordan chain of $L$ at $\lambda=0$.
Namely, observing that $Z:=\partial_\lambda U$ satisfies at $\lambda=0$
the first-order generalized eigenvalue equation
\be\label{gen}
(L-\lambda)Z= U
\ee
if $U(\lambda, \cdot)$ satisfies for all $\lambda$ the eigenvalue
equation $(L-\lambda) U=0$,
we may arrange that not only 
$W^+_1(0)= W^-_4(0)=\partial_x \bar W^{\bar s}$ as above, but also
\be\label{chain}
\partial_\lambda W^+_1(0)= \partial_\lambda W^-_4(0)=
-(\partial \bar W^s/\partial s)(\bar s).
\ee
For, $W_1^+$ and $W_4^-$ may be chosen as ``fast modes'', decaying
uniformly exponentially as $x\to \pm \infty$ for $|\lambda|$ sufficiently
small, whence $\partial_\lambda W_1^+(0)$ and $\partial_\lambda W_4^-(0)$
are uniformly exponentially decaying at their associated spatial infinities
as well, and (their $U$ components) 
satisfy the generalized eigenvalue equation \eqref{gen}, which
uniquely specifies them up to exponentially decaying solutions of
the eigenvalue equation: in this case, multiples  $c_1W_1^+$ and $c_2W_4^-$,
respectively, which may be removed by the analytic change of coordinates
$W_j\to \lambda c_j W_j$.

With this normalization, we find immediately that 
$$
D(0)=
\det \bp \partial_x \bar W^{\bar s}  & W_2^+& W_3^-
&\partial_x \bar W^{\bar s}  \ep =0
$$
and
\ba
D'(0)&=
\det \bp \partial_\lambda W^+_1  & W_2^+& W_3^-
& W_4^-  \ep +\cdots\\
&\quad +
\det \bp W^+_1  & W_2^+& W_3^-
& \partial_\lambda W_4^-  \ep \\
&=
\det \bp \partial_x \bar W^{\bar s}  & W_2^+& W_3^-
&\partial_s \bar W^{\bar s}  \ep 
\\
&\quad +
\det \bp \partial_s \bar W^{\bar s}  & W_2^+& W_3^-
&\partial_x \bar W^{\bar s}  \ep 
 =0.
\ea

By a similar computation, we find that
\ba\label{D''}
D''(0)&=
\det \bp \partial_\lambda^2 W^+_1  & W_2^+& W_3^-
& W_4^-  \ep \\
&\quad +
\det \bp W^+_1  & W_2^+& W_3^-
& \partial_\lambda^2 W_4^-  \ep\\
&=
\det \bp 
W_1^+& W_2^+& W_3^-
& 
(\partial_\lambda^2 W_4^-  
- \partial_\lambda^2 W_1^+  )
\ep,\\
\ea
$W_1^+= \partial_x \bar W^{\bar s} $,
where $Y_j:=\partial_\lambda^2 U_j$
satisfy at $\lambda=0$ the second-order
generalized eigenvalue equation
\be\label{2gen}
(L-\lambda)Y= Z.
\ee
(Recall that $Z:=\partial_\lambda U$ satisfies the first-order generalized
eigenvalue equation \eqref{gen}.)

Now, setting $\lambda=0$ and integrating \eqref{eval} from $x=\pm \infty$
to $x=0$, we obtain
\ba\label{inteval}
sv + u&=(sv+u)(\pm \infty), \\
-su + \alpha v + \kappa v''&= (-su + \alpha v)(\pm \infty),\\
\ea
for each $(u,v)=(u, v)^\pm_j$ associated with $W_j^\pm$.
In particular, the righhand sides of \eqref{inteval} vanish for
$(u,v)=(u,v)_1^+$.
Likewise, setting 
$(\tilde u, \tilde v)=
(\tilde u, \tilde v)_j^\pm :=
 \partial_\lambda^2(u, v)^\pm_j$,
$j=1, 4$, using \eqref{2gen}, and recalling that we chose $(\tilde u, \tilde v)$
to vanish at spatial infinity,
we find that 
\ba\label{intvar}
s\tilde v + \tilde u&= -\int_{\pm \infty}^x
(\partial \bar v^{\bar s}/\partial s)\, dx,\\
-s\tilde u + \alpha \tilde v + \kappa \tilde v''&= 
\int_{\pm \infty}^x
(\partial \bar u^{\bar s}/\partial s)\, dx.\\
\ea

Performing the row operations corresponding the the lefthand side
of \eqref{inteval}, i.e., adding $\kappa^{-1}(-s\tilde u + \alpha \tilde v)$ 
to $v''$, then adding $sv$ to $u$, and using \eqref{inteval},
\eqref{lims}, and \eqref{intvar}, reduces \eqref{D''} to
\ba\label{redD''}
D''(0)&=
\det \bp 
0 & s+c& s-c& C_v\\
 v_1^+& v_2^+& v_3^- & 
(\partial_\lambda^2 v_4^-  
- \partial_\lambda^2 v_1^+  )\\
 (v_1^+)'& (v_2^+)'& (v_3^-)' & 
(\partial_\lambda^2 v_4^-  
- \partial_\lambda^2 v_1^+  )'\\
0 & -sc-c^2& sc-c^2& C_u\\
\ep,\\
\ea
where
$C_v:=-\int_{- \infty}^{+\infty} (\partial \bar v^{\bar s}/\partial s)\, dx$
and
$C_u:=\kappa^{-1}
\int_{- \infty}^{+\infty} (\partial \bar u^{\bar s}/\partial s)\, dx$.
Since
$$
\bp
s+c& s-c\\
-sc-c^2& sc-c^2\\
\ep
$$
is invertible, by \eqref{signalpha}, there exist constants
$d_2$, $d_3$ such that
$$
d_3
\bp
s-c\\
sc-c^2\\
\ep
-d_2
\bp
s+c\\
-sc-c^2\\
\ep
=
\bp
C_v\\
C_u\\
\ep.
$$
Performing the corresponding column operation to eliminate
$C_v$, $C_u$, we obtain, finally,
\ba\label{fredD''}
D''(0)&=
\det \bp 
0 & s+c& s-c& 0\\
 v_1^+& v_2^+& v_3^- & 
(\partial_\lambda^2 v_4^-  
- \partial_\lambda^2 v_1^+ + d_2 v_2^+-d_3 v_3^-  )\\
 (v_1^+)'& (v_2^+)'& (v_3^-)' & 
(\partial_\lambda^2 v_4^-  
- \partial_\lambda^2 v_1^+ + d_2 v_2^+-d_3 v_3^-  )'\\
0 & -sc-c^2& sc-c^2& 0\\
\ep,\\
&=
C\gamma, 
\ea
where 
\ba\label{C}
C=\det \bp
s+c& s-c\\
-sc-c^2& sc-c^2\\
\ep
 %=2c(-c^2+s^2)
=2c(p'(v_\infty)+s^2)<0.
\ea
and
\be
\gamma=
\det \bp 
 v_1^+ & 
(\partial_\lambda^2 v_4^-  
 \partial_\lambda^2 v_1^+ + d_2 v_2^+-d_3 v_3^-  )\\
 (v_1^+)'  & 
%(\hat v_4-\hat v_1)'\\
(\partial_\lambda^2 v_4^-  
 \partial_\lambda^2 v_1^+ + d_2 v_2^+-d_3 v_3^-  )'\\
\ep.\\
\ee

Expand now $\gamma=\gamma_-(0)-\gamma_+(0)$, where
\be\label{gammas}
\gamma_-(x):=
\det \bp 
 v_1^+ & 
\hat v_4\\
 (v_1^+)' & 
(\hat v_4)'\\
\ep,
\qquad
\gamma_+(x):=
\det \bp 
 v_1^+ & 
\hat v_1\\
 (v_1^+)' & 
(\hat v_1)'\\
\ep\\
\ee
and 
\be\label{hatv}
\hat v_4:= \partial_\lambda^2 v_4^+ +  d_3 v_3^-  ,
\qquad
 \hat v_1:=  \partial_\lambda^2 v_1^+ + d_2 v_2^+.
\ee
Since $\hat v_-$ is bounded as $x\to -\infty$, while $v_4^-$
decays exponentially, 
\be\label{vanish1}
\gamma_-(-\infty)=0
\ee
Similarly,
\be\label{vanish2}
\gamma_+(+\infty)=0.
\ee

By \eqref{inteval}, \eqref{intvar}, $v_1^+=\bar v_x$ satisfies
$$
 v'' + \kappa^{-1}(s^2  + \alpha )v=0,
$$
or
\be\label{fv}
\bp v \\  v'\ep'-
\bp
0 & 1 \\
-\kappa^{-1}(s^2  + \alpha )& 0\\
\ep
\bp v \\  v'\ep
=\bp 0\\0 \ep
\ee
while $\hat v_\pm$ satisfy
$$
\hat v'' + \kappa^{-1}(s^2  + \alpha )\hat v=
\kappa^{-1}\int_{- \infty}^x 
[-s (\partial \bar v^{\bar s}/\partial s)
+ (\partial \bar u^{\bar s}/\partial s)]\, dx
+\hat C,
$$
or
\be\label{fhatv}
\bp \hat v \\  \hat v'\ep'-
\bp
0 & 1 \\
-\kappa^{-1}(s^2  + \alpha )& 0\\
\ep
\bp \hat v \\  \hat v'\ep
=\bp 0\\ 
F
 \ep,
\ee
with $C$ constant,
$$
F:=
\kappa^{-1}\int_{- \infty}^x 
[-s (\partial \bar v^{\bar s}/\partial s)
+ (\partial \bar u^{\bar s}/\partial s)]\, dx
+\hat C.
$$

Thus, using inhomogeneous Abel's formula and 
\eqref{vanish1}-\eqref{vanish2}, we may evaluate the difference
$\gamma$ between Wronskians
$\gamma_+$ and $\gamma_-$ at $x=0$ as the Melnikov-type integral
\ba\label{Mel}
\gamma&=\gamma_-(0)-\gamma_+(0)\\
&=
\int_{-\infty}^{+\infty}
\det
\bp v_1^+ & 0\\
0 & F
\ep (x) \, dx \\
&=
\kappa^{-1}\int_{-\infty}^{+\infty}
\bar v_x(x)
\int_{- \infty}^x 
[-s (\partial \bar v^{\bar s}/\partial s)
+ (\partial \bar u^{\bar s}/\partial s)](y)\, dy \, dx\\
&=
\kappa^{-1}\int_{-\infty}^{+\infty}
(\bar v- v_\infty)
[-s (\partial \bar v^{\bar s}/\partial s)
+ (\partial \bar u^{\bar s}/\partial s)](x)\, dx\\
\ea
where in the second-to-last line we are using the fact that term
$\hat C \int_{-\infty}^{+\infty} v_x(x)\, dx$ integrates to zero
and in the last line integration by parts.

Substituting the relation
$-s(\bar v-v_\infty)=(\bar u-u_\infty)$
coming from the first equation in \eqref{trav} and recalling
\eqref{monotone}, Remark \ref{mon}, we obtain, finally,
\ba\label{gform}
\gamma&=
\kappa^{-1}\int_{-\infty}^{+\infty}
[(\bar u-u_\infty)
(\partial \bar v^{\bar s}/\partial s)
+
(\bar v -v _\infty)
(\partial \bar u^{\bar s}/\partial s)](x)\, dx\\
&=
\kappa^{-1}\int_{-\infty}^{+\infty}
(\bar U-U_\infty)J
(\partial \bar U^{\bar s}/\partial s)(x)\, dx\\
&=
\kappa^{-1}
(d^2m/ds^2)(\bar s).
\ea
Combining \eqref{fredD''} and \eqref{gform}, we obtain
\be\label{signD''}
D''(0)=(-C/\kappa) (d^2m/ds^2)(\bar s),
\ee
$-C/\kappa > 0$,
completing the proof. 
\end{proof}

\subsection{Proof of Corollary \ref{stab}}

By standard considerations \cite{PW, GZ, BSZ, Z2} we have also the following.

\begin{lem}\label{infsign}
As $\lambda \to +\infty$ along the real axis,
$\sgn D(\lambda)$ has limit $+1$.
\end{lem}

\begin{proof}
First, note that
$D(\lambda)$ does not vanish for $\Re \lambda$ sufficiently
large, a standard fact associated with well-posedness 
of the linearized time-evolutionary system;
this may be established either by asymptotic ODE theory as in
\cite{GZ, Z2} or by elementary energy estimates as in \cite{BSZ}.
Thus, $D(\lambda)$ has a (nonzero) limiting sign as claimed.

To determine the value of this sign, we may consider a homotopy
from system \eqref{eval} to the constant-coefficient equation
\ba\label{homotopy}
\lambda v  - u'&=0\\
\lambda u  + \kappa v'''&=0,\\
\ea
capturing high-frequency behavior, for which bases 

$$
V_j=\bp \lambda/\mu_j, 1, \mu_j, \mu_j^2 \ep,
\qquad \mu_j=
\theta_j \sqrt \lambda 
$$ 
of exponential
solutions may be explicitly calculated for all $\lambda$,
where $\theta_j$ are the fourth roots of $-1$.

Choosing bases
$V_1$, $(V_2-V_1)/\sqrt \lambda$ and
$V_3$, $ (V_4-V_3)/\sqrt \lambda$ of decaying subspaces
at $\pm \infty$, we find that these 
extend continuously to $\lambda=0$,
with the projection onto $(u,v)$ coordinates of the limiting bases equal 
to the standard Euclidean basis $(1,0)$, $(0,1)$; in particular, the projection
of the limiting subspaces is nonsingular.
Likewise, it is nonsingular for the bases at $\lambda=0$ for the original
system, and (by identical calculation) for all systems strictly between the two.
Thus, choosing continuously initializing bases at $\lambda =0$ for
the family of systems, we find that the signs of $D(+\infty)$
and of the determinants $d_\pm$ of the projections of decaying subspaces
at $\pm \infty$ onto $(u,v)$ coordinates are all constant under homotopy, 
with, moreover $\sgn D(+\infty)d_- d_+$ explicitly calculable from
system \eqref{homotopy}.
We may therefore determine $\sgn D(+\infty)$ for the original system
by the straightforward calculation of $d_-$ and $d_+$, yielding the result.
%We omit the details.
% for similar arguments, see \cite{Z2}.
%TODO:
%Actually, this is more general I think!  Can we do this way general
%MP stable viscosity for n\times n???? (!)
%A- no! As we lack the result at $\lambda=0$... ???
\end{proof}

\begin{rem}\label{alt}
\textup{
An alternative, somewhat simpler approach, using the homoclinic 
structure of the wave, is to note that the limiting
subspaces at plus and minus infinity of decaying solutions of the
eigenvalue equation are complementary subspaces of the limiting
coefficient matrix $A$, hence transverse.  Thus, their determinant
(real-valued, by construction) is of constant sign independent of $\lambda$,
which may be related to $\sgn D(+\infty)$ by explicit calculation.
Calculating the value at $\lambda=0$ of this determinant then gives the result.
This is similar in spirit to the original argument of \cite{PW}.
However, this approach does not extend to the heteroclinic case, and
so we have chosen to give the more general argument above.
For related arguments, see \cite{Z2, Z3}.
}
\end{rem}

We may now easily obtain the main result.

\begin{proof}[Proof of Corollary \ref{stab}]
It has been shown by variational considerations in \cite{BDDJ}
that $(d^2m/ds^2)(\bar s)>0$ implies linearized and nonlinear stability,
using Corollary \ref{bound} and the argument of \cite{GSS}.
Thus, we need only establish that
$(d^2m/ds^2)(\bar s)<0$ implies instability, for which it is sufficient
to show that there is an $L^2$ eigenfunction of $L$ with positive real
eigenvalue $\lambda$.
Defining the {\it stability index}
\be\label{Gamma}
\Gamma:= \sgn D''(0) D(+\infty),
\ee
we have, provided $\Gamma\ne 0$, that the number of positive
real eigenvalues of $L$ (zeroes of $D$) is odd or even, according
as $\Gamma$ is negative or positive,
by standard degree
theory on the line (recall, $D$ is real-valued, by construction,
for $\lambda$ real).
Thus, $\Gamma<0$ implies existence of at least one positive real eigenvalue,
for which (see Section \ref{s:evans}) 
the associated eigenfunction is necessarily exponentially decaying
as $x\to \pm \infty$.
%Combining these facts, we find that the number of unstable eigenvalues
%is one if $\Gamma<0$ and zero if $\Gamma >0$.
Noting that $\sgn \Gamma= \sgn (d^2m/ds^2)(\bar s)$
by \eqref{signD''} combined with Lemma \ref{infsign},
we are done.
%by the proof of Theorem \ref{main} (which includes the additional
%information that $\sgn c=\sgn D(+\infty)$), we are done.
%by Theorem \ref{main}, we find that $\bar U$ is spectrally stable
%if $(d^2m/ds^2)(\bar s)>0$ and only if $(d^2m/ds^2)(\bar s)\ge 0$; 
%in particular, it is linearly and nonlinearly {\it unstable} if
%$(d^2m/ds^2)(\bar s)<0$. 
%Recalling the result already shown by variational considerations in \cite{BDDJ}
%that $(d^2m/ds^2)(\bar s)>0$ implies linearized and nonlinear stability,
%we are done. 
\end{proof}

\begin{rem}\label{count}
\textup{
By Lemma \ref{PWbound}, 
there is in fact precisely one unstable eigenvalue,
necessarily real, in case $(d^2m/ds^2)(\bar s)<0$. 
}
\end{rem}

\section{Remarks on the viscous case}\label{s:viscous}

Finally, we comment briefly on the viscous case, as
modeled by
\ba\label{viscous}
v_t - u_x&=0\\
u_t + p(v)_x&= \epsilon u_{xx}
-\kappa v_{xxx},\\
\ea
$\epsilon >0$.  
Stability of traveling waves of this model has been studied in detail
in \cite{Z, OZ} for a parameter range of $(\kappa, \epsilon)$ for
which \eqref{viscous} may be converted by a change of dependent
variables to a strictly parabolic problem.

%\begin{rem}
%\textup{
%The physical significance for dynamical phase transitions is less clear.
Including viscosity of arbitrarily small strength
$\epsilon>0$, one finds by the energy argument of \cite{OZ}
that all homoclinic connections of the $\epsilon=0$ equations (2) break, 
save for zero-speed solitary waves identical to those of \eqref{2}.
Moreover, by a calculation similar to but much simpler than that
of the previous section, we find that all of these waves are unstable,
with $D(0)=0$, $D'(0)=c \epsilon \int (\bar v_x)^2 dx $,
$c\ne 0$.  That is, even infinitesimal viscosity will either break or
destabilize any solitary wave, whether stable or unstable for the
$\epsilon=0$ model, leaving unclear the physical implications of stability
or instability with respect to \eqref{2}.   
However, depending whether the $\epsilon=0$
version is stable or unstable, the unstable root will be near or
far from the origin, corresponding perhaps to some type of metastable
phenomenon.
%}
%\end{rem}

\bigbreak

\end{document}